\documentstyle[12pt]{article}
\begin{document}
\baselineskip=16pt

\centerline{}
\vskip0.5cm
\centerline{\Large\bf Convex Hypersurfaces and}
\vskip0.4cm
\centerline{\Large\bf $L^p$ Estimates for Schr\"odinger Equations}

\vskip0.8cm
\centerline{\sc Quan Zheng}
\centerline{Department of Mathematics}
\centerline{Huazhong University of Science and Technology}
\centerline{Wuhan 430074, P. R. China}
\centerline{and}
\centerline{Center for Optimal Control and Discrete Mathematics}
\centerline{Huazhong Normal University}
\centerline{Wuhan 430079, P. R. China}
\centerline{qzheng@hust.edu.cn}

\vskip0.4cm
\centerline{Xiaohua Yao}
\centerline{Department of Mathematics}
\centerline{Huazhong University of Science and Technology}
\centerline{Wuhan 430074, P. R. China}
\centerline{yaoxiaohua@hust.edu.cn}

\vskip0.4cm
\centerline{\sc Dashan Fan}
\centerline{Department of Mathematics}
\centerline{University of Wisconsin-Milwaukee}
\centerline{Milwaukee, WI 53201, USA}
\centerline{fan@csd.uwm.edu}

\renewcommand {\thefootnote}{} \footnote{\small

\vspace{-0.4cm}\noindent
This project was supported by the National Science
Foundation of China
}

\rightskip1cm\leftskip1cm
{\small
\vskip0.5cm
\centerline{\bf Abstract}

\vskip0.2cm
This paper is concerned with Schr\"odinger equations whose
principal operators are homogeneous elliptic. When the
corresponding level hypersurface is convex, we show the
$L^p$-$L^q$ estimate of solution operator in free case. This
estimate, combining with the results of fractionally integrated
groups, allows us to further obtain the $L^p$ estimate of
solutions for the initial data belonging to a dense subset of
$L^p$ in the case of integrable potentials. }

\rightskip0cm\leftskip0cm

\newpage
\baselineskip=18pt

\centerline{}
\vskip1cm
\centerline{1. \ {\sc Introduction}}

\vskip0.4cm
In this paper we take interest in $L^p$, $1\le p<\infty$,
estimates of solutions for the following Schr\"odinger equation
$$
{\partial u\over \partial t}=(iP(D)+V)u,\quad u(0,\cdot)=u_0\in
L^p({\bf R}^n),
\eqno(*)$$
where $D=-i(\partial/\partial x_1\cdots,\partial/\partial x_n)$,
$P:{\bf R}^n\to{\bf R}$ is a homogeneous elliptic polynomial of
order $m$ ($m$ must be even, except $n=1$), and $V$ is a suitable
potential function. In the sequel, we may assume without loss of
generality that $P(\xi)>0$ for $\xi\ne0$. Otherwise we have
$P(\xi)<0$ for $\xi\ne0$, for which the following hypersurface
$\Sigma$ should be replaced by
$$
\{\xi\in{\bf R}^n|\ P(\xi)=-1\}.
$$

In order to obtain $L^p$ estimates of the solution of $(*)$, we
will first treat $L^p$-$L^q$ estimates of $e^{itP(D)}$, which is
the solution operator of $(*)$ with $V=0$. To this end, we need to
consider the compact hypersurface
$$
\Sigma=\{\xi\in{\bf R}^n|\ P(\xi)=1\}.
$$
When the Gaussian curvature of $\Sigma$ is nonzero everywhere, it
is known that $L^p$-$L^q$ estimates of $e^{itP(D)}$ $(t\ne0)$ can
be deduced from Miyachi [12]. In fact, Miyachi gave some remarks
on these estimates in a more general case where $P$ is a positive
and smooth homogeneous function, provided the nonvanishing
Gaussian curvature on $\Sigma$. Also, dropping the homogeneity of
$P$, Balabane and Emami-Rad [4] studied these estimates under a
suitable nondegenerate condition. However, one can check that the
nondegenerate condition is equivalent to the nonzero Gaussian
curvature if $P$ is homogeneous.

As we know, the nonvanishing Gaussian curvature plays a crucial
rule in estimating many oscillatory integrals [16]. This is a
reason why one needs such condition in [12, 4]. Howerve, there
exist many hypersurfaces $\Sigma$ whose Gaussian curvatures may
vanish at some points (although we have observed that if $m=2$
then $\Sigma$ has nonzero Gaussian curvature everywhere under our
assumptions on $P$). These examples are easily available, for
instance, the hypersurfaces $\Sigma$ associated with polynomials
$\xi_1^m+\cdots+\xi_n^m$ $(m=4,6,\cdots)$ or $\xi_1^4+6\xi_1^2
\xi_2^2+\xi_2^4$.

On the other hand, an important subclass of hypersurfaces with
vanishing Gaussian curvature at some points is the class of convex
hypersurfaces of finite type [5]. The main purpose of this paper
is to investigate the $L^p$ estimate of the solution of $(*)$ when
$\Sigma$ is a convex hypersurface of finite type. Roughly
speaking, this means that $P$ allows to be degenerate on a subset
of ${\bf R}^n$.

This paper is organized as follows.

In section 2, we study $L^p$-$L^q$ estimates of the solution
operator $e^{itP(D)}$ $(t\ne0)$ and the resolvent operator
$(\lambda -iP(D))^{-1}$ (Re$\lambda\ne0$) when $\Sigma$ is a
convex hypersurface of finite type. The method used is quite
different from those in the previous papers [12, 4], due to the
nature of the vanishing Gaussian curvature. Our proof depends
heavily on a decay estimate for the kernels ${\cal F}^{-1}(e^{\pm
iP})$, in which we need to use a powerful theorem in [5]. Since
the proof is involved and very technical, we will present it in
section 3.

In section 4, we show that the operator $iP(D)+V$ with suitable
integrable potential $V$ generates an integrated group on
$L^p({\bf R}^n)$. As we know, the semigroup of operators is a
useful abstract tool to treat Cauchy problems. However, the Cauchy
problem $(*)$ in $L^p({\bf R}^n)$ $(p\ne2)$ cannot be treated by
classical semigroups of operators (i.e., $C_0$-semigroups). In
fact, the Schr\"odinger operator $iP(D)$ generates a
$C_0$-semigroups in $L^p({\bf R}^n)$ if and only if $p=2$ (see
[10, 12]). Thus, several generalizations of $C_0$-semigroups, such
as smooth distribution semigroups [3], integrated semigroups [1,
9], and regularized semigroups [7, 8] were introduced and applied
to different general differential operators [9, 18]. In our case,
We use fractionally integrated groups to deal with the Cauchy
problem $(*)$ in $L^p({\bf R}^n)$, which will leads to better
results than using smooth distribution semigroups (see [3]).
Moreover, when $P$ is nondegenerate type $m$, we will show how our
results present an improvement over Theorem 2$'$ and Theorem 6 in
[4].

Throughout this paper, denote by $\Sigma$ the hypersurface
$\{\xi\in{\bf R}^n|\ P(\xi)=1\}$. Assume, except in the last
section, that $P:{\bf R}^n\to[0,\infty)$ is always a homogeneous
elliptic polynomial of order $m$ where $n\ge2$, $m$ is even and
$\ge4$.

\vskip1cm
\centerline{2. \ {\sc $L^p$-$L^q$ estimates for Schr\"odinger
equations without potentials}}

\vskip0.4cm

We start with the concept of finite type. Denote by $S$ the smooth
hypersurface $\{\xi\in{\bf R}^n|\ \phi(\xi)=0\}$, where $\phi\in
C^\infty({\bf R}^n)$ and $\nabla \phi(\xi)\ne0$ for $\xi\in S$. We
say that $S$ is of finite type if any one dimensional tangent line
has at most a finite order of contact with $S$. The precise
definition is as follows.

Denote by ${\bf S}^{n-1}$ the unit sphere in ${\bf R}^n$. Let
$$
\nabla_\eta=\sum_{j=1}^n\eta_j\partial/\partial x_j \quad{\rm
for}\ \eta=(\eta_1,\cdots,\eta_n)\in{\bf S}^{n-1},
$$
which is the directional derivative in direction $\eta$, and let
$\nabla_\eta^j$ be the $j$-th power of this derivative.

\vskip0.4cm
{\bf Definition 2.1.} \ {\it Let $k$ be an integer. The smooth
hypersurface $S$ is of type $k$ if there exists a constant
$\delta>0$ such that
$$
\sum_{j=1}^k|\nabla_\eta^j\phi(\xi)|\ge\delta\quad {\it for}\
\xi\in S\ {\it and}\ \eta\in{\bf S}^{n-1}.
$$
Moreover, we say that $S$ is convex if
$$
S\subset\{\eta\in{\bf R}^n|\ \langle\eta-\xi,\nabla\phi(\xi)
\rangle\ge0\}\quad {\it for}\ \xi\in S
$$
or}
$$
S\subset\{\eta\in{\bf R}^n|\ \langle\eta-\xi,\nabla\phi(\xi)
\rangle\le0\}\quad {\it for}\ \xi\in S.
$$

\vskip0.4cm
It is clear that $k\ge2$, and that if $S$ is of type $k$ it is
also of type $k'(>k)$. For the hypersurface $\Sigma$ (i.e.
$\{\xi\in{\bf R}^n|\ P(\xi)=1\}$), since
$$
\langle\xi,\nabla P(\xi)\rangle=mP(\xi)=m\quad{\rm for}\
\xi\in\Sigma,
$$
it follows that $\nabla P(\xi)\ne0$ for $\xi\in\Sigma$, and thus
$\Sigma$ is smooth. Also, a simple computation leads to
$$
\nabla_\eta^m(P(\xi)-1)=m!P(\eta)\quad{\rm for}\ \xi\in\Sigma\
{\rm and}\ \eta\in{\bf S}^{n-1}.
$$
Hence we have

\vskip0.4cm
{\bf Proposition 2.2.} \ {\it $\Sigma$ is a smooth compact
hypersurface of type less and equal to $m$.}

\vskip0.4cm
A simple example of polynomials whose level hypersurface $\Sigma$
is of type $m$ is $\xi_1^m+\cdots+\xi_n^m$ $(m=4,6,\cdots)$. We
notice that there exist polynomials $P$ whose level hypersurfaces
$\Sigma$ are of type $k(<m)$. For example, when $P(\xi)=\xi_1^6+
5\xi_1^2\xi_2^4+\xi_2^6$ the corresponding hypersurface $\Sigma$
is of type 4, but $m=6$.

We now turn to the Cauchy problem $(*)$ with $V=0$. In this case,
for every initial data $u_0\in{\cal S}({\bf R}^n)$ (the Schwartz
space), the solution is given by
$$
u(t,\cdot)=e^{itP(D)}u_0:={\cal F}^{-1}(e^{itP})*u_0,
$$
where ${\cal F}$ (or $\hat{}$\,) denotes the Fourier transform,
${\cal F}^{-1}$ its inverse, and ${\cal F}^{-1}(e^{itP})$ is
understood in the distributional sense. Therefore, to obtain
$L^p$-$L^q$ estimates of $e^{itP(D)}$ $(t\ne0)$, the key result is
to show estimates of the kernels ${\cal F}^{-1}(e^{\pm iP})$.

In the sequel, denote by $p'$ the conjugate index of $p$, and
$\|\cdot\|_{L^p\mbox{-}L^q}$ the norm in ${\cal L}(L^p,L^q)$ (the
space of all bounded linear operators from $L^p$ to $L^q$). Let
$$
h(m,n,k)={m-2\over2(m-1)}+{(m-k)(n-1)\over k(m-1)}\quad{\rm for}\
2\le k\le m,
$$
$\tau=n/h(m,n,k)$, and $q(p)=q(m,n,k,p)$ where
$$
{1\over q(m,n,k,p)}={1\over\tau p}+{1\over\tau'p'}\quad{\rm for}\
1\le p<2.
$$
We first remark that when $2\le k\le m$,
$$
{2(m-1)\over m-2}\le\tau\le{2n(m-1)\over m-2}.
$$
Since $m\ge4$, it follows that $\tau\in(2,3n]$. Next, we remark
that
$$
{1\over p}>{1\over\tau p}+{1\over\tau'p'}>{1\over p'}\quad{\rm
for}\ 1\le p<2,
$$
and thus $2<q(p)<p'$. Moreover, denote by $I_p$ ($1\le p\le2$) the
following subset of $[2,\infty]$:
$$
I_p=\cases{(q(p),\infty] &${\rm if}\ 1\le p<\tau',$\cr
(q(p),{p(2-\tau')\over p-\tau'}) &${\rm if}\ \tau'\le p<2,$\cr
\{2\} &${\rm if}\ p=2.$}
$$

\vskip0.4cm
{\bf Theorem 2.3.} \ {\it Suppose $\Sigma$ is a convex
hypersurface of type $k$. Then ${\cal F}^{-1}(e^{\pm iP})\in
C^\infty({\bf R}^n)$ and}
$$
({\cal F}^{-1}(e^{\pm iP}))(x)=O(|x|^{-h(m,n,k)})\quad{\it as}\
|x|\to\infty.
$$

\vskip0.4cm The proof is lengthy and is given in the next section.

\vskip0.4cm
{\bf Theorem 2.4.} \ {\it Suppose $\Sigma$ is a convex
hypersurface of type $k$. If $p\in[1,2]$ and $q\in I_p$, then
there exists a constant $C>0$ such that}
$$
\|e^{itP(D)}\|_{L^p\mbox{-}L^q}\le C|t|^{{n\over m}({1\over q}
-{1\over p})}\quad{\it for}\ t\ne0.
$$

\vskip0.4cm
{\bf Proof.} \ By Theorem 2.3, ${\cal F}^{-1}(e^{\pm iP})\in
L^s({\bf R}^n)$ for $s>\tau$. Since $P$ is homogeneous, one has
$$
{\cal F}^{-1}(e^{itP})(x)=|t|^{-n/m}{\cal F}^{-1}(e^{itP/|t|})
(|t|^{-1/m}x)\quad{\rm for}\ t\ne0\ {\rm and}\ x\in{\bf R}^n,
$$
and thus
$$
\|{\cal F}^{-1}(e^{itP})\|_{L^s}=|t|^{-n/ms'}\|{\cal F}^{-1}
(e^{itP/|t|})\|_{L^s}\le C|t|^{-n/ms'}\quad{\rm for}\ t\ne0,
$$
where the constant $C$ is independent of $t$. The remaining of the
proof will be divided into several steps.

Step 1. When $1\le p<\tau'$ and ${\tau'p\over\tau'-p}<q\le\infty$,
it follows from Young's inequality that
$$
\|e^{itP(D)}\|_{L^p\mbox{-}L^q}\le\|{\cal F}^{-1}(e^{itP})\|_{L^s}
\le C|t|^{{n\over m}({1\over q}-{1\over p})}\quad{\it for}\ t\ne0,
$$
where ${1\over s}=1+{1\over q}-{1\over p}$, which implies
$s>\tau$.

Step 2. Since $P(D)$ is selfadjoint in $L^2({\bf R}^n)$,
$\|e^{itP(D)}\|_{L^2\mbox{-}L^2}=1$ for $t\ge0$ by Stone's
theorem. When $1\le p<2$ and $q(p)<q\le p'$, we deduce from the
Riesz-Thorin interpolation theorem that
$$
\|e^{itP(D)}\|_{L^p\mbox{-}L^q}\le\|e^{itP(D)}\|_{L^1\mbox{-}L^s}
^{1-2/p'}\|e^{itP(D)}\|_{L^2\mbox{-}L^2}^{2/p'}\le C|t|^{{n\over
m}({1\over q}-{1\over p})}\quad{\it for}\ t\ne0,
$$
where $s={q(p'-2)\over p'-q}>\tau$.

Step 3. When $1\le p<\tau'$ and $q(p)<q\le\infty$, we notice
$q(p)<{\tau'p\over\tau'-p}$. Since $L^{s_2}({\bf R}^n)\subset
L^{s_1}({\bf R}^n)+L^{s_3}({\bf R}^n)$ for $1\le s_1\le s_2\le
s_3\le\infty$, the desired estimate is a direct consequence of the
conclusions in Steps 1 and 2.

Step 4. When $\tau'\le p<2$ and $q(p)<q<{p(2-\tau')\over
p-\tau'}$, we put $\lambda={2(p-p_0)\over p(2-p_0)}$ where
$p_0\in[1,\tau')$ such that $q<{2\over\lambda}<{p(2-\tau')\over
p-\tau'}$. A simple computation leads to
$$
{\lambda\over2}<{1\over q}<{1\over q(p)}={1-\lambda\over q(p_0)}
+{\lambda\over2}.
$$
Consequently, there exists a unique $q_0>q(p_0)$ such that
${1\over q}={1-\lambda\over q_0}+{\lambda\over2}$. Also, ${1\over
p}={1-\lambda\over p_0}+{\lambda\over2}$. The desired estimate now
can be deduced from the Riesz-Thorin interpolation theorem and the
conclusion in Step 1. $\hfill\Box$

\vskip0.4cm
The subsequent theorem deals with $L^p$-$L^q$ estimates of the
resolvent of $iP(D)$.

\vskip0.4cm
{\bf Theorem 2.5.} \ {\it Suppose $\Sigma$ is a convex
hypersurface of type $k$. If $p\in[1,2]$, $q\in I_p$, and ${1\over
p}-{1\over q}<{m\over n}$, then there exists a constant $C>0$ such
that}
$$
\|(\lambda-iP(D))^{-1}\|_{L^p\mbox{-}L^q}\le C|{\rm Re}\lambda|
^{{n\over m}({1\over p}-{1\over q})-1}\quad{\it for}\ {\rm Re}
\lambda\ne0.
$$

\vskip0.4cm
{\bf Proof.} \ For Re$\lambda>0$ and $f\in{\cal S}({\bf R}^n)$,
one has
\begin{eqnarray*}
(\lambda-iP(D))^{-1}f&=&{\cal F}^{-1}((\lambda-iP)^{-1}\hat{f})\\
&=&\int_0^\infty e^{-\lambda t}{\cal F}^{-1}(e^{itP}\hat{f})dt\\
&=&\int_0^\infty e^{-\lambda t}e^{itP(D)}fdt.
\end{eqnarray*}
It follows therefore from Theorem 2.4 that
\begin{eqnarray*}
\|(\lambda-iP(D))^{-1}\|_{L^p\mbox{-}L^q}
&\le&C\int_0^\infty e^{-({\rm Re}\lambda)t}t^{{n\over m}
({1\over q}-{1\over p})}dt\\
&=&C|{\rm Re}\lambda|^{{n\over m}({1\over p}-{1\over q})-1}.
\end{eqnarray*}

For Re$\lambda<0$ and $f\in{\cal S}({\bf R}^n)$, one has
$$
(\lambda-iP(D))^{-1}f=-(-\lambda+iP(D))^{-1}f=\int_0^\infty
e^{\lambda t}e^{-itP(D)}fdt,
$$
and thus the desired estimate follows from Theorem 2.4.
$\hfill\Box$

\vskip0.4cm
Since $p'\in I_p$ for $p\in[1,2]$, we have

\vskip0.4cm
{\bf Corollary 2.6.} \ {\it Suppose $\Sigma$ is a convex
hypersurface of type $k$. If $p\in[1,2]$ then
$$
\|e^{itP(D)}\|_{L^p\mbox{-}L^{p'}}\le C|t|^{{n\over m}(1-{2\over
p})}\quad{\it for}\ t\ne0.
$$
If, in addition, $p>{2n\over n+m}$ then}
$$
\|(\lambda-iP(D))^{-1}\|_{L^p\mbox{-}L^{p'}}\le C|{\rm Re}\lambda|
^{-{2n\over mp'}}\quad{\it for}\ {\rm Re} \lambda\ne0.
$$

\vskip1cm
\centerline{3. \ {\sc The proof of Theorem 2.3}}

\vskip0.4cm
By our assumptions on $P$, $\phi:=P^{1/m}$ is a positive and
smooth homogeneous function of degree 1, and $\Sigma=\{\xi\in{\bf
R}^n|\ \phi(\xi)=1\}$. Let $\varphi\in C^\infty({\bf R})$ such
that supp$\varphi\subset[1,\infty)$ and $\varphi(t)=1$ for $t\ge
2$. Obviously, in order to estimate ${\cal F}^{-1}(e^{iP})$
(similarly for ${\cal F}^{-1}(e^{-iP})$), it suffices to estimate
${\cal F}^{-1}((\varphi\circ\phi)e^{iP})$ (cf. [4, p.363]).
Consider the integral
\begin{eqnarray*}
K_\varepsilon(x)&:=&\int_{{\bf R}^n}e^{-\varepsilon\phi(y)+iP(y)
+i\langle x,y\rangle}\varphi(\phi(y))dy \cr
&=&\int_0^\infty e^{-\varepsilon t+it^m}t^{n-1}\varphi(t)
\Big{(}\int_\Sigma{e^{itr\langle\eta,\xi\rangle}\over|\nabla
\phi(\xi)|}d\sigma(\xi)\Big{)}dt\quad{\rm for}\ \varepsilon>0,
\end{eqnarray*}
where $r=|x|$, $x=r\eta$, and $d\sigma$ is the induced surface
measure on $\Sigma$. We will show that
$$
{\cal F}^{-1}((\varphi\circ\phi)e^{iP})(x)=(2\pi)^{-n}
\lim_{\varepsilon \to0}K_\varepsilon(x)
$$
uniformly for $x$ in compact subsets of ${\bf R}^n$, and
$K_\varepsilon(x)$ decays as $|x|^{-h(m,n,k)}$. From this we have
${\cal F}^{-1}((\varphi\circ\phi)e^{iP})\in C({\bf R}^n)$ since it
is clear that $K_\varepsilon\in C({\bf R}^n)$.

Denote by $\Pi$ the Gaussian map
$$
\xi\in\Sigma\longmapsto{\nabla\phi(\xi)\over|\nabla\phi(\xi)|}
\in{\bf S}^{n-1}.
$$
Since $\Sigma$ is a compact convex hypersurface, $\Pi$ is a
homeomorphism from $\Sigma$ to ${\bf S}^{n-1}$. For given
$\eta\in{\bf S}^{n-1}$, let $\xi_\pm=\Pi^{-1}(\pm\eta)$. Then
$$
\langle\eta,\xi_\pm\rangle=\pm\Big{\langle}{\nabla\phi(\xi_\pm)
\over|\nabla\phi(\xi_\pm)|},\xi_\pm\Big{\rangle}=\pm{\phi(\xi_\pm)
\over|\nabla\phi(\xi_\pm)|}=\pm{1\over|\nabla\phi(\xi_\pm)|},
$$
Noting that $\pm\eta$ is the outward unit normal to $\Sigma$ at
$\xi_\pm$, by Theorem B in [5] (also cf. [6]) we have
$$
\int_\Sigma{e^{i\lambda\langle\eta,\xi\rangle}\over|\nabla
\phi(\xi)|}d\sigma(\xi)=e^{i\lambda\langle\eta,\xi_+\rangle}
H_+(\lambda)+e^{i\lambda\langle\eta,\xi_-\rangle}H_-(\lambda)
+H_\infty(\lambda)\quad{\rm for}\ \lambda>0.
$$
Here $H_{\pm}\in C^\infty((0,\infty))$,
$$
|H_{\pm}^{(j)}(\lambda)|\le C_j\lambda^{-j-(n-1)/k}\quad{\rm for}\
j\in{\bf N}_0,
$$
and
$$
|H_\infty(\lambda)|\le C_j\lambda^{-j}\quad{\rm for}\ j\in{\bf N},
$$
where ${\bf N}_0={\bf N}\cup\{0\}$ and constants $C_j$ depend only
on the hypersurface $\Sigma$. Hence
\begin{eqnarray*}
K_\varepsilon(x)&=&\int_0^\infty e^{-\varepsilon t+it^m
+itr\langle\eta,\xi_+\rangle}t^{n-1}\varphi(t)H_+(tr)dt \cr
&&+\int_0^\infty e^{-\varepsilon t+it^m+itr\langle\eta,
\xi_-\rangle}t^{n-1}\varphi(t)H_-(tr)dt \cr
&&+\int_0^\infty e^{-\varepsilon t+it^m}t^{n-1}\varphi(t)
H_\infty(tr)dt \cr
&:=&J_1^\varepsilon+J_2^\varepsilon+J_3^\varepsilon.
\end{eqnarray*}
In the remaining of this section, for the sake of convenience, we
will denote by $C$ a generic constant independent of $r$, $t$ and
$\varepsilon$.

We consider first the integral $J_3^\varepsilon$. It is obvious
that
$$
|J_3^\varepsilon|\le C\int_1^\infty t^{n-1}(tr)^{-(n+1)}dt\le
Cr^{-(n+1)}.
$$
Since $r=|x|$, it follows from the dominated convergence theorem
that $J_3^\varepsilon$ ($\varepsilon\to0$) converges uniformly for
$x$ in compact subsets of ${\bf R}^n\setminus\{0\}$, and decays as
$|x|^{-h(m,n,k)}$, where we notice that $h(m,n,k)\le n$.

Next, consider the integral $J_1^\varepsilon$. Let
$$
\cases{u(t)=-\varepsilon t+it^m+itr\langle\eta,\xi_+\rangle\cr
v(t)=t^{n-1}\varphi(t)H_+(tr)^{^{^{^{}}}}}
$$
for $t>0$. Since $u'(t)\ne0$ for $t>0$, we can define
$D_{\mbox{\tiny \#}}f=-gf'$ and $D_*f=(gf)'$ for $f\in C^1
((0,\infty))$, where $g=-1/u'$. By induction on $j$ we find that
$$
g^{(j)}(t)=\sum_{l=l_0}^ja_lt^{l(m-1)-j}g(t)^{l+1}\quad{\rm for}\
j\in{\bf N}_0,
$$
where constants $a_l$ depend only on $l$ and $m$, and $l_0\in{\bf
N}_0$ such that $l_0\ge j/(m-1)$. Since there exists a constant
$c>1$ such that $c^{-1}\le|\nabla\phi(\xi)|\le c$ for $\xi\in
\Sigma$, we have
$$
|g(t)|=1/|u'(t)|\le1/r\langle\eta,\xi_+\rangle=|\nabla\phi(\xi_+)|/r
\le c/r\quad{\rm for}\ t>0.
$$
Also, $|g(t)|\le{1\over m}t^{1-m}$ for $t>0$. Hence
$$
|g^{(j)}(t)|\le Cr^{-1}t^{-j}\quad{\rm for}\ j\in {\bf N}_0.
$$
On the other hand, one sees
$$
\Big{|}{d^j\over dt^j}(H_+(tr))\Big{|}\le Ct^{-j}(tr)^{-(n-1)/k}
\quad{\rm for}\ j\in{\bf N}_0,
$$
and thus by Leibniz's formula
$$
|v^{(j)}(t)|\le Cr^{-(n-1)/k}t^{-j+n-1-(n-1)/k}\quad{\rm for}\
j\in {\bf N}_0.
\eqno(3.1)$$
Since it is not hard to show
$$
D_*^jv=\sum_\alpha a_\alpha g^{(\alpha_1)}\cdots g^{(\alpha_j)}
v^{(\alpha_{j+1})}\quad{\rm for}\ j\in {\bf N},
$$
where the sum runs over all $\alpha=(\alpha_1,\cdots,\alpha_{j+1})
\in{\bf N}_0^{j+1}$ such that $|\alpha|=j$ and $0\le\alpha_1\le
\cdots\le\alpha_j$, it follows that
$$
|(D_*^jv)(t)|\le Cr^{-j-(n-1)/k}t^{-j+n-1-(n-1)/k}\quad{\rm for}\
j\in {\bf N}_0, \eqno(3.2)$$ where we used the notation
$D_*^0v=v$. Noting that $D_{\mbox{\tiny \#}}^ne^u=e^u$ we have
$$
J_1^\varepsilon=\int_0^\infty(D_{\mbox{\tiny \#}}^ne^u)(t)v(t)dt
=\int_0^\infty e^{u(t)}(D_*^nv)(t)dt.
$$
Consequently
$$
|J_1^\varepsilon|\le Cr^{-n-(n-1)/k}\int_1^\infty
t^{-1-(n-1)/k}dt\le Cr^{-n-(n-1)/k}.
$$
The dominated convergence theorem yields thus that
$J_1^\varepsilon$ converges uniformly for $x$ in compact subsets
of ${\bf R}^n\setminus\{0\}$, and decays as $|x|^{-h(m,n,k)}$.

We now consider the integral $J_2^\varepsilon$. In this case, we
put
$$
\cases{u(t)=-\varepsilon t+it^m-it\bar{r}\cr
v(t)=t^{n-1}\varphi(t)H_-(tr)^{^{^{^{}}}}}
$$
for $t>0$, where
$$
\bar{r}:=-r\langle\eta,\xi_-\rangle=r/|\nabla\phi(\xi_-)|.
$$
Since $t_0:=(\overline{r}/m)^{1/(m-1)}$ is the unique critical
point of the oscillatory integral $J_2^\varepsilon$, we write
\begin{eqnarray*}
J_2^\varepsilon&=&\Big{\{}\int_{2t_0}^\infty+\int_{t_0/2}^{2t_0}
+\int_0^{t_0/2}\Big{\}}e^{u(t)}v(t)dt\\
&:=&J_{2,1}^\varepsilon+J_{2,2}^\varepsilon+J_{2,3}^\varepsilon.
\end{eqnarray*}

From integration by parts one gets
$$
J_{2,1}^\varepsilon=-{e^{u(2t_0)}\over u'(2t_0)}\sum_{j=0}^{n-1}
(D_*^jv)(2t_0)+\int_{2t_0}^\infty e^{u(t)}(D_*^nv)(t)dt.
$$
Since
$$
|u'(t)|\ge mt^{m-1}-\bar{r}\ge(2^{m-1}-1)\bar{r}\ge(2^{m-1}-1)r/c
\quad{\rm for}\ t\ge2t_0
$$
and since
$$
|u'(t)|\ge mt^{m-1}-\bar{r}\ge m(1-2^{1-m})t^{m-1}\quad{\rm for}\
t\ge2t_0,
$$
the estimate (3.2) still holds for $t\ge2t_0$. Hence
\begin{eqnarray*}
|J_{2,1}^\varepsilon|&\le&Cr^{-1}\sum_{j=0}^{n-1}r^{-j-(n-1)/k}
(2t_0)^{-j+n-1-(n-1)/k}\\
&&+C\int_{2t_0}^\infty r^{-n-(n-1)/k}t^{-1-(n-1)/k}dt\\
&\le&Cr^{(n-m-m(n-1)/k)/(m-1)}\sum_{j=0}^{n-1}r^{-jm/(m-1)}.
\end{eqnarray*}
But
$$
(n-m-m(n-1)/k)/(m-1)\le-h(m,n,k),
$$
$J_{2,1}^\varepsilon$ converges uniformly for $x$ in compact
subsets of ${\bf R}^n\setminus\{0\}$, and decays as
$|x|^{-h(m,n,k)}$. Since $|u'(t)|\ge Ct^{m-1}$ and $\ge Cr$ for
$0<t\le t_0/2$, a slight modification of the above method leads to
the same conclusion for $J_{2,3}^\varepsilon$. We omit the
details.

In order to deal with $J_{2,2}^\varepsilon$, set
$w(t)=t^m-t\bar{r}$ and $v(t)$ is defined as in
$J_{2,1}^\varepsilon$. Obviously,
$$
\lim_{\varepsilon\to0}J_{2,2}^\varepsilon=J_{2,2}:=\int_{t_0/2}
^{2t_0}e^{iw(t)}v(t)dt
$$
uniformly for $x$ in compact subsets of ${\bf R}^n\setminus\{0\}$.
It remains to show that $J_{2,2}$ decays as $|x|^{-h(m,n,k)}$. Let
$y=(t-t_0)/t_0$, which maps $[t_0/2,2t_0]$ onto $[-1/2,1]$, and
let $\lambda=m(m-1)t_0^m$. Then
\begin{eqnarray*}
\Phi(y)&:=&{1\over\lambda}(w(t_0(y+1))-w(t_0))\\
&=&{1\over m(m-1)}((y+1)^m-my-1)\\
&=&{1\over m(m-1)}\sum_{l=2}^m{m\choose l}y^l\quad{\rm for}\
y\in[-1/2,1].
\end{eqnarray*}
Consequently
$$
\Phi''(y)=(y+1)^{m-2}\ge2^{2-m}\quad{\rm for}\ y\in[-1/2,1]
$$
and
$$
J_{2,2}=t_0e^{iw(t_0)}\int_{-1/2}^1e^{i\lambda\Phi(y)}v(t_0(y+1))dy.
$$
It follows thus from van der Corput's theorem (cf. [16, p.334])
that
\begin{eqnarray*}
|J_{2,2}|&\le&C\lambda^{-1/2}t_0\Big{(}|v(2t_0)|+\int_{-1/2}^1
|t_0v'(t_0(y+1))|dy\Big{)}\\
&\le&Ct_0^{1-m/2}\Big{(}|v(2t_0)|+\sup_{t_0/2\le t\le2t_0}
|t_0v'(t)|\Big{)}.
\end{eqnarray*}
Since (3.1), in which $H_+$ is replaced by $H_-$, still holds,
$$
|J_{2,2}|\le Cr^{-(n-1)/k}t_0^{n-m/2-(n-1)/k}\le Cr^{-h(m,n,k)},
$$
as desired.

Finally, we only need to show that $K_\varepsilon(x)$ converges
uniformly for $x$ in some neighborhood of the origin. Let $U$ be
the ball $\{x\in{\bf R}^n|\ |x|\le m/2M\}$, where $M$ is a
constant such that $|\xi|\le M$ for $\xi\in\Sigma$. For given
$x\in U$, let
$$
u(t)=-\varepsilon t+it^m+it\langle x,\xi\rangle\quad{\rm for}\
t\ge1.
$$
Then
$$
K_\varepsilon(x)=\int_\Sigma\Big{(}\int_1^\infty e^{u(t)}t^{n-1}
\varphi(t)dt\Big{)}{d\sigma(\xi)\over|\nabla\phi(\xi)|}.
$$
Since
$$
|u'(t)|\ge|mt^{m-1}+\langle x,\xi\rangle|\ge mt^{m-1}-M|x|
\ge(m/2)t^{m-1}\ge m/2,
$$
the same way as in the estimate of $J_1^\varepsilon$ yields
$$
\Big{|}\int_1^\infty e^{u(t)}t^{n-1}\varphi(t)dt\Big{|}\le C,
$$
and therefore the claim follows.

Moreover, an analogous method as above leads to
$$
D^\alpha{\cal F}^{-1}(e^{\pm iP})={\cal F}^{-1}(\xi^\alpha e^{\pm
iP(\xi)})\in C({\bf R}^n)\quad{\rm for}\ \alpha\in{\bf N}^n_0,
$$
i.e. ${\cal F}^{-1} (e^{\pm iP})\in C^\infty({\bf R}^n)$ and
$$
D^\alpha({\cal F}^{-1}(e^{\pm iP}))(x)=O(|x|^{-h(m,n,k)
+|\alpha|/(m-1)})\ (|x|\to\infty)\quad{\rm for}\ \alpha\in{\bf
N}^n_0.
$$

\vskip1cm
\centerline{4. \ {\sc $L^p$ estimates for Schr\"odinger
equations}}

\vskip0.4cm
It was Balabane and Emami-Rad [4] who first applied smooth
distribution semigroups to higher order Schr\"odinger equations,
and showed the $L^p$ estimate of solutions. Arendt and Kellermann
[2] showed that a smooth distribution semigroup is equivalent to
an integrally integrated semigroup. However, it is known that the
fractionally integrated semigroup is a generalization of the
integrally integrated semigroup and is a more suitable tool for
elliptic differential operators in $L^p({\bf R}^n)$ (cf. [9, 17]).
We start with its definition.

\vskip0.4cm
{\bf Definition 4.1.} \ {\it Let $A$ be a linear operator on a
Banach space $X$ and $\beta\ge 0$. Then a strongly continuous
family $T:[0,\infty)\to{\cal L}(X)$ is called a $\beta$-times
integrated semigroup on $X$ with generator $A$ if there exist
constants $C,\omega\ge0$ such that $\|T(t)\|\le Ce^{\omega t}$ for
$t\ge 0$, $(\omega,\infty)\subset\rho(A)$ (the resolvent set of
$A$), and
$$
(\lambda-A)^{-1}x=\lambda^\beta\int_0^\infty e^{-\lambda t}T(t)x
dt\quad{\it for}\ \lambda>\omega\ {\it and}\ x\in X.
$$

If $A$ and $-A$ both are generators of $\beta$-times integrated
semigroups on $X$, we say $A$ is the generator of a $\beta$-times
integrated group on $X$.}

\vskip0.4cm
Here is a sufficient condition for an operator to be the generator
of integrated semigroups, which is due to Hieber [9, p.30] (see
[2] for a special case).

\vskip0.4cm
{\bf Lemma 4.2.} \ {\it Let $A$ be a linear operator on a Banach
space $X$. Suppose there exist constants $C,\omega\ge0$ and
$\gamma\ge-1$ such that for ${\rm Re}\lambda>\omega$, $\lambda
\in\rho(A)$ and $\|(\lambda-A)^{-1}\|\le C|\lambda|^\gamma$. Then
$A$ generates a $\beta$-times integrated semigroup on $X$, where
$\beta>\gamma+1$.}

\vskip0.4cm
Assume that the operator $P(D)$ has maximal domain in the
distributional sense on $L^p({\bf R}^n)$ ($1\le p<\infty$), and
thus it is closed and densely defined. It is known that $D(P(D))W^{m,p}({\bf R}^n)$ for $1<p<\infty$. The following Lemma 4.3(a)
can be found in [9, 17], and Lemma 4.3(b) follows from Lemma
4.3(a) and Definition 4.1, immediately.

\vskip0.4cm
{\bf Lemma 4.3.} \ {\it Let $1\le p<\infty$ and $\beta>n_p:n|{1\over2}-{1\over p}|$.

(a) $iP(D)$ generates a $\beta$-times integrated group $T(t)$
$(t\in {\bf R})$ on $L^p({\bf R}^n)$, and
$$
\|T(t)\|_{L^p\mbox{-}L^p}\le C|t|^\beta\quad{\it for}\ t\in{\bf
R}.
$$

(b) $\{\lambda\in{\bf C}|\ {\rm Re}\lambda\ne0\}\subset
\rho(iP(D))$ on $L^p({\bf R}^n)$, and}
$$
\|(\lambda-iP(D))^{-1}\|_{L^p\mbox{-}L^p}\le C|\lambda|^\beta/
|{\rm Re}\lambda|^{\beta+1}\quad{\it for}\ {\rm Re}\lambda\ne0.
$$

\vskip0.4cm
Let $V$ be a measurable function defined on ${\bf R}^n$. We
consider $V$ as a multiplication operator on $L^p({\bf R}^n)$ with
$D(V):=\{f\in L^p({\bf R}^n)|\ Vf\in L^p({\bf R}^n)\}$. The domain
of $iP(D)+V$ is $D(P(D))\cap D(V)$. Denote by $I'_p$ ($1\le
p\le2$) the following subset of $[1,\infty]$:
$$
I'_p=\cases{[p,{\tau'p\over 2-p}) &${\rm if}\ 1\le p<\tau',$\cr
({p(2-\tau')\over 2-p},{\tau'p\over 2-p}) &${\rm if}\ \tau'\le
p<2,$\cr \{\infty\} &${\rm if}\ p=2,$}
$$
where $\tau=n/h(m,n,k)$ and $h(m,n,k)$ is defined as in section 2.

\vskip0.4cm
{\bf Theorem 4.4.} \ {\it Let $V=V_1+V_2$ with $V_j\in L^{s_j}
({\bf R}^n)$ for some $s_j\in({n\over m},\infty]$ $(j=1,2)$.

(a) If $s_j\in I'_p$ for some $p\in[1,2]$, then $iP(D)+V$
generates a $\beta$-times integrated group on $L^p({\bf R}^n)$,
where $\beta>n_p+1$.

(b) If $s_j\in I'_{p'}$ for some $p\in(2,\infty)$, then an
extension of $iP(D)+V$, i.e. $(-iP(D)+\overline{V})^*$ generates a
$\beta$-times integrated group on $L^p({\bf R}^n)$, where
$\beta>n_p+1$.}

\vskip0.4cm
{\bf Proof.} \ Since $iP(D)+V$ and $-(iP(D)+V)$ satisfy the same
assumptions, it suffices to show that $iP(D)+V$ generates a
$\beta$-times integrated semigroup on $L^p({\bf R}^n)$.

We consider first the case $1\le p\le2$. Let ${1\over q_j}={1\over
p}-{1\over s_j}$ ($j=1,2$). Then $s_j\in I'_p$ implies $q_j\in
I_p$ ($I_p$ is defined in section 2). So we obtain by Theorem 2.5
and H\"older's inequality that
\begin{eqnarray*}
\|V(\lambda-iP(D))^{-1}\|_{L^p\mbox{-}L^p}
&\le&\sum_{j=1,2}\|V\|_{L^{q_j}\mbox{-}L^p}\|(\lambda-iP(D))^{-1}\|
_{L^p\mbox{-}L^{q_j}}\\
&\le& C\sum_{j=1,2}\|V\|_{L^{s_j}}|{\rm Re}\lambda|^{{n\over
ms_j}-1}.
\end{eqnarray*}
In view of ${n\over ms_j}-1<0$, there exists $\omega\ge1$ such
that
$$
\|V(\lambda-iP(D))^{-1}\|_{L^p\mbox{-}L^p}\le1/2\quad{\rm for\
Re}\lambda>\omega.
$$
Consequently, $\lambda\in\rho(iP(D)+V)$ and
$$
(\lambda-iP(D)-V)^{-1}=(\lambda-iP(D))^{-1}\sum_{j=0}^\infty
(V(\lambda-iP(D))^{-1})^j.
$$
This implies by Lemma 4.3(b) that
\begin{eqnarray*}
\|(\lambda-iP(D)-V)^{-1}\|_{L^p\mbox{-}L^p}
&\le&2\|(\lambda-iP(D))^{-1}\|_{L^p\mbox{-}L^p} \\
&\le&C|\lambda|^{n_p+\varepsilon}\quad{\rm for\ Re}\lambda>\omega,
\end{eqnarray*}
where $\varepsilon\in(0,\beta-n_p-1)$. It follows now from Lemma
4.2 that $iP(D)+V$ generates a $\beta$-times integrated semigroup
on $L^p({\bf R}^n)$ .

Next, we consider the case $2<p<\infty$. From the proof of (a) one
sees that $-iP(D)+\overline{V}$ is densely defined on $L^{p'}({\bf
R}^n)$, and thus $(-iP(D)+\overline{V})^*$ exists and is densely
defined on $L^p({\bf R}^n)$. It is easy to check $iP(D)+V\subset
(-iP(D)+\overline{V}) ^*$. Also, an adjointness argument implies
$$
\|(\lambda-(-iP(D)+\overline{V})^*)^{-1}\|_{L^p\mbox{-}L^p}
=\|(\overline{\lambda}-(-iP(D)+\overline{V}))^{-1}\|_{L^{p'}
\mbox{-}L^{p'}}.
$$
Since $s_j\in I'_{p'}$ and $n_{p'}=n_p$, this leads to the same
estimate as in the case $1<p<2$. It follows therefore from Lemma
4.2 that $(-iP(D)+ \overline{V})^*$ generates a $\beta$-times
integrated semigroup on $L^p({\bf R}^n)$. $\hfill\Box$

\vskip0.4cm
When $2<p<\infty$, we rewrite $I'_{p'}$ as
$$
I'_{p'}=\cases{({p(2-\tau')\over p-2},{\tau'p\over p-2}) &${\rm
if}\ 2<p\le\tau,$\cr [{p\over p-1},{\tau'p\over p-2}) &${\rm if}\
\tau<p< \infty.$}
$$
From this and the subsequent proposition one sees that it is not
always true that $iP(D)+V=(-iP(D)+\overline{V})^*$ in Theorem
4.3(b), so is the such situation in [4, Theorem 6]. This means
that the operator $iP(D)+V$ in [4, Theorem 6] should be replaced
by $(-iP(D)+\overline{V})^*$ for $p>2$.

\vskip0.4cm
{\bf Proposition 4.5.} \ {\it Let $1<p<\infty$ and ${n\over m}<s<
\infty$.

(a) If $V\in L^s({\bf R}^n)$ and $s\ge p$, then $W^{m,p}({\bf
R}^n)\subset D(V)$ and $V$ is a compact operator from
$W^{m,p}({\bf R}^n)$ to $L^p({\bf R}^n)$.

(b) If $V\in L^s({\bf R}^n)$ and $s\ge\max\{p,p'\}$, then
$iP(D)+V=(-iP(D)+\overline{V})^*$ on $L^p({\bf R}^n)$.

(c) If $1\le s<p$, then there exists $V\in L^s({\bf R}^n)$ such
that $W^{m,p}({\bf R}^n)\cap D(V)=\{0\}$.}

\vskip0.4cm
{\bf Proof.} \ (a) is a direct consequence of Theorem 10.2 in [15,
p.147], in which the condition (5) is satisfied. Since $\rho(P(D))
\ne\emptyset$ on $L^{p'}({\bf R}^n)$, (b) follows from (a) and
Theorem 6.1 in [15, p.94]. To show (c), we will modify slightly an
example in [14, p.60]. Let $f(x)=|x|^{-n/p}$ for $|x|\le1$ and
$=0$ for $|x|>1$. Define
$$
V(x)=\sum_{j=1}^\infty2^{-j}f(x-\alpha_j)\quad{\rm for}\ x\in{\bf
R}^n,
$$
where $\{\alpha_j\}_{j=1}^\infty={\bf Q}^n$ (${\bf Q}$ denotes the
set of all rational numbers). Then
$$
\|V\|_{L^s}\le\sum_{j=1}^\infty2^{-j}\|f(x-\alpha_j)\|_{L^s}
=\|f\|_{L^s}<\infty.
$$
If $0\ne g\in W^{m,p}({\bf R}^n)\cap D(V)$, one has $g\in C({\bf
R}^n)$ by Sobolev's embedding theorem. Thus $|g(x)|\ge c>0$ in
some open subset $\Omega\subset{\bf R}^n$. Taking $\alpha_j\in
\Omega$ yields
$$
\int_{{\bf R}^n}|V(x)g(x)|^pdx
\ge2^{-np}c^p\int_{|x-\alpha_j|\le\delta}|x-\alpha_j|^{-n}dx
=\infty\quad{\rm for\ small}\ \delta>0,
$$
which contradicts $g\in D(V)$. $\hfill\Box$

\vskip0.4cm
Set $I'_p=I'_{p'}\cap[p,\infty)$ for $p\in(2,\infty)$. Then
$I'_p=\emptyset$ for $p\ge2+\tau'$ and
$$
I'_p=\cases{({p(2-\tau')\over p-2}, {\tau'p\over p-2}) &${\rm if}\
2<p\le4-\tau',$\cr [p,{\tau'p\over p-2}) &${\rm if}\
4-\tau'<p<2+\tau'.$}
$$
Combining Theorem 4.4 and Proposition 4.5(b) leads to

\vskip0.4cm
{\bf Theorem 4.6.} \ {\it Let $1\le p<2+\tau'$ and $V=V_1+V_2$
with $V_j\in L^{s_j}({\bf R}^n)$ for some $s_j\in I'_p\cap({n\over
m},\infty]$ $(j=1,2)$. Then $iP(D)+V$ generates a $\beta$-times
integrated group on $L^p({\bf R}^n)$, where $\beta>n_p+1$.}

\vskip0.4cm
Corresponding to Corollary 2.6 we have

\vskip0.4cm
{\bf Corollary 4.7.} \ {\it Let $1\le p\le3$, $n_p<m/2$ and $V\in
L^{p\over|p-2|}({\bf R}^n)$. Then $iP(D)+V$ generates a
$\beta$-times integrated group on $L^p({\bf R}^n)$, where
$\beta>n_p+1$.}

\vskip0.4cm
In order to give $L^p$ estimates of the solution for Schr\"odinger
equations we need Straub's fractional powers (cf. [13]). If a
densely defined operator $A$ is the generator of a $\beta$-times
integrated group $T(t)$ $(t\in{\bf R})$ satisfying $\|T(t)\|\le
Ce^{\omega|t|}$ $(t\in{\bf R})$. Then for $\delta,\varepsilon>0$
the fractional powers $(\omega+\delta\pm A)^{\beta+\varepsilon}$
are well-defined and their domains are independent of $\sigma>0$.
We note that $D((\omega+\delta+A)^{\beta+\varepsilon})\cap
D((\omega+\delta-A)^{\beta+\varepsilon})$ for small
$\varepsilon>0$ contains the dense subspace $D(A^{[\beta]+1})$.
The following result is a consequence of Theorem 4.6 and Theorem
1.1 in [13].

\vskip0.4cm
{\bf Theorem 4.8.} \ {\it Suppose $p$, $V$ and $\beta$ satisfy the
assumptions of Theorem 4.6. Then there exist constants
$C,\omega>0$ such that for every data $u_0\in D((\omega+iP(D)+V)
^\beta)\cap D((\omega-iP(D)-V)^\beta)$, the Cauchy problem $(*)$
has a unique solution $u\in C({\bf R},L^p({\bf R}^n))$ and
$$
\|u(t,\cdot)\|_{L^p}\le Ce^{\omega|t|}\|(\omega\pm iP(D)\pm V)
^\beta u_0\|_{L^p}\quad{\it for}\ t\in{\bf R},
$$
where we choose $+\ ($resp. $-)$ if $t\ge0\ ($resp. $<0)$.}

\vskip0.4cm
When $P$ is nondegenerate (i.e. det$\Big{(}{\partial^2P(\xi)\over
\partial\xi_i\partial\xi_j}\Big{)}_{n\times n}\ne0$ for $\xi\in
{\bf R}^n\setminus\{0\}$), the Gaussian curvature of $\Sigma$ is
nonzero everywhere (cf. [11]). In this case $k=2$, and thus
$h(m,n,k)= {n(m-2)\over2(m-1)}$. In order to compare our results
with those in [4] for homogeneous polynomial $P$, we denote by
$I_p(\tau)$ (resp. $I'_p(\tau)$) the set $I_p$ (resp. $I'_p$)
defined in section 2 (resp. 4). As a direct consequence of
Theorems 2.4 and 4.4(a) we have

\vskip0.4cm
{\bf Corollary 4.9.} \ {\it Suppose $P$ is nondegenerate. Let
$1\le p\le2$ and $\tau_0={2(m-1)\over m-2}$. If $q\in I_p(\tau_0)$
$($resp. $s_j\in I'_p(\tau_0)\cap({n\over m},\infty])$, then the
conclusion of Theorem 2.4 $($resp. 4.4$($a$))$ is true.}

\vskip0.4cm
We first note that for homogeneous polynomial $P$, the hypothesis
(H2) in [4] is equivalent to that $P$ is nondegenerate. Thus
Corollary 4.9 improves the corresponding Theorems 2$'$ and 6 in
[4] in several respects:

1) For fixed $p\in[1,2)$, the interval $I_p(\tau_0)$ is replaced
by smaller $(q(\tau_1,p),p']$ in Theorem 2$'$, where ${1\over
q(\tau,p)}={1\over\tau p}+{1\over\tau'p'}$ and $\tau_1={2n(m-1)
\over mn-2n-3m+2}$. In fact, it is clear that $q(\tau,p)$ $(\tau>
2)$ is strictly increasing. Since $p'\in I_p(\tau_0)$ and $2<
\tau_0<\tau_1$, $I_p(\tau_0)$ contains properly $(q(\tau_1,p),
p']$. Similarly, the interval $I'_p(\tau_0)$ is replaced by
smaller $[{p\over 2-p},{\tau'_1p\over 2-p})$ in Theorem 6.

2) The hypothesis (H3$'$), i.e. $n>3+{4\over m-2}$ is required in
Theorems 2$'$ and 6, but not in Corollary 4.9. For example, when
$n=2,3$ and $m\ge4$ one has $n\le3+{4\over m-2}$, and thus
Theorems 2$'$ and 6 can not deal with such case. However, in this
case, $I'_p(\tau_0)\cap({n\over m},\infty]=I'_p(\tau_0)$, which
means that for every $p\in[1,2]$, we can choose $q$ and $s_j$'s
values such that the conclusion of Corollary 4.9 holds.

3) We first note that $p\ne1$ in Theorems 2$'$ and 6, but it is
admitted that $p=1$ in Corollary 4.9. Furthermore, it is required
in Theorem 6 that $p>{2c\over c+1}$, where $c$ is an integer with
$c>{n\over m-1}$. However, the restriction of $p$ in Corollary 4.9
is only caused by $I'_p(\tau_0)\cap({n\over m},\infty]\ne
\emptyset$, which is equivalent to $p>{2n\over n+2m-2}$. It is
easy to see that this is naturally an improvement of the
corresponding condition in Theorem 6.

4) The conclusion in Theorem 6 (see its remark for homogeneous
$P$) is that $iP(D)+V$ generates a smooth distribution group on
$L^p({\bf R}^n)$ of order $\beta$, which is equivalent to a
$\beta$-times integrated group on $L^p({\bf R}^n)$, where $\beta$
is an integer with $\beta>n_p+2$. Our conclusion in Corollary 4.9
however admits that $\beta$ is a real number with $\beta>n_p+1$.


\vskip1cm
\begin {thebibliography}{20}
\bibitem{1} W. Arendt, {\it Vector-valued Laplace transforms and Cauchy problems},
       Israel J. Math. {\bf 59} (1987) 327-352.
\bibitem{2} W. Arendt and H. Kellermann, {\it Integrated solutions of Volterra
       integro-differential equations and applications}, In: Volterra Integro-differential
       Equations in Banach Spaces and Applications (G. Da Prato and M. Iannelli, eds.),
       Longman, Harlow, 1989, 21-51.
\bibitem{3} M. Balabane and H. A. Emami-Rad, {\it Smooth distribution group and
       Schr\"odinger equation in $L^p$}, J. Math. Anal. Appl. {\bf 70} (1979) 61-71.
\bibitem{4} M. Balabane and H. A. Emami-Rad, {\it $L^p$ estimates for Schr\"odinger
       evolution equations}, Trans. Amer. Math. Soc. {\bf 292} (1985) 357-373.
\bibitem{5} J. Bruna, A. Nagel, and S. Wainger, {\it Convex hypersurfaces and
       Fourier transforms}, Ann. Math. {\bf 127} (1988), 333-365.
\bibitem{6} M. Cowling, S. Disney, G. Maucceri, and D. M\"uller,
       {\it Damping oscillatory integrals}, Invent. Math. {\bf 101} (1990) 237-260.
\bibitem{7} G. Da Prato, {\it Semigruppi regolarizzabili}, Ricerche Mat. {\bf 15}
       (1966) 223-248.
\bibitem{8} E. B. Davies and M. M. H. Pang, {\it The Cauchy problem and a
       generalization of the Hille-Yosida theorem}, Proc. London Math. Soc.
       {\bf 55} (1987) 181-208.
\bibitem{9} M. Hieber, {\it Integrated semigroups and differential operators
       on $L^p$}, Dissertation, T\"ubingen, 1989.
\bibitem{10} L. H\"ormander, {\it Estimates for translation invariant operators in
       $L^p$ spaces}, Acta Math. {\bf 104} (1960) 93-140.
\bibitem{11} A. Miyachi, {\it On some estimates for the wave equation in $L^p$ and $H^p$},
       J. Fac. Sci. Univ. Tokyo {\bf 27} (1980) 331-354.
\bibitem{12} A. Miyachi, {\it On some singular Fourier multipliers}, J. Fac. Sci.
       Univ. Tokyo {\bf 28} (1981) 267-315.
\bibitem{13} J. van Neerven and B. Straub, {\it On the existence and growth of mild
       solutions of the abstract Cauchy problem for operators with polynomially bounded
       resolvent}, Houston J. Math. {\bf 24} (1998) 137-171.
\bibitem{14} M. Schechter, {\it Operator Methods in Quantum Mechanics},
       Elsevier North Holland, New York, 1981.
\bibitem{15} M. Schechter, {\it Spectra of Partial Differential Operators},
       2nd ed., Elsevier Science Publishers B.V., Amsterdam, 1986.
\bibitem{16} E. M. Stein, {\it Harmonic Analysis: Real-Variable Methods, Orthogonality,
       and Oscillatory Integrals}, Princeton Univ. Press, New Jersey, 1993.
\bibitem{17} Q. Zheng, {\it Abstract differential operators and Cauchy problems},
       T\"ubinger Berichte zur Funktionalanalysis {\bf 4} (1995) 273-280.
\bibitem{18} Q. Zheng and Y. Li, {\it Abstract parabolic systems and regularized
       semigroups}, Pacific J. Math. {\bf 182} (1998) 183-199.
\end {thebibliography}
\end{document}